\documentclass[a4paper,12pt,reqno]{amsart}
\usepackage{latexsym}
\usepackage{amssymb} 
\usepackage{mathrsfs}
\usepackage{amsmath}
\usepackage{latexsym}
\usepackage{delarray}
\usepackage{amssymb,amsmath,amsfonts,amsthm,mathrsfs}
\usepackage{blkarray}

\usepackage{hyperref}
\renewcommand\eqref[1]{(\ref{#1})} 
 \newtheorem{thm}{Theorem}[section]
 \newtheorem{cor}[thm]{Corollary}
 \newtheorem{lem}[thm]{Lemma}
 \newtheorem{prop}[thm]{Proposition}
 \theoremstyle{definition}
 \newtheorem{defn}[thm]{Definition}
 \theoremstyle{remark}

 \numberwithin{equation}{section}
\newcommand{\half}{\frac{1}{2}}

\newcommand{\ene}{\mathbb{N}}

\newcommand{\er}{\mathbb{R}}
\newcommand{\ce}{\mathbb{C}}

\newcommand{\efee}{\mathcal{F}}

\newcommand{\bi}{\begin{itemize}}
\newcommand{\cinfm}{{C}^{\infty}(M)}

\newcommand{\ei}{\end{itemize}}
\newcommand{\be}{\begin{enumerate}}
\newcommand{\ee}{\end{enumerate}}
\newcommand{\beq}{\begin{equation}}
\newcommand{\eq}{\end{equation}}

\newcommand{\Dcal}{\mathcal{D}}

\newcommand{\subll}{\mathcal{L}_{sub}}
\def\SO3{{{\rm SO(3)}}}
\newcommand{\lapm}{\Delta_M}

\DeclareMathOperator{\Tr}{Tr}

\def\N{{{\mathbb N}}}

\def\SU2{{{\rm SU(2)}}}
\def\lapsu2{{{\mathcal L}_\SU2}}

\begin{document}

\title[Kernel and symbol criteria for Schatten classes and $r$-nuclearity]
 {Kernel and symbol criteria for Schatten classes and $r$-nuclearity on compact manifolds}

\author[Julio Delgado]{Julio Delgado}

\address{%
Department of Mathematics\\
Imperial College London\\
180 Queen's Gate, London SW7 2AZ\\
United Kingdom
}
\email{j.delgado@imperial.ac.uk}

\thanks{The first author was supported by Marie Curie IIF 301599. 
The second author was supported by EPSRC grant EP/K039407/1. }
\author[Michael Ruzhansky]{Michael Ruzhansky}

\address{%
Department of Mathematics\\
Imperial College London\\
180 Queen's Gate, London SW7 2AZ\\
United Kingdom
}

\email{m.ruzhansky@imperial.ac.uk}

\date{\today}
\begin{abstract}
In this note we present criteria on both symbols and integral kernels ensuring
 that the corresponding operators on compact manifolds belong to Schatten classes. 
A specific test for nuclearity is established as well as the corresponding trace formulae. 
In the special case of compact Lie groups, kernel criteria in terms of (locally and globally)
 hypoelliptic operators are also given. A notion of an invariant operator and 
 its full symbol associated to an elliptic operator are introduced. 
 Some applications to the study of $r$-nuclearity on $L^{p}$ spaces are also obtained.
\end{abstract}

\maketitle

\section{Introduction}

Let $M$ be a closed smooth manifold (smooth manifold without boundary)
endowed with a positive measure $dx$. We denote by $\Psi^{\nu}(M)$ the 
usual H\"ormander class of pseudo-differential 
 operators of order $\nu\in\er$.
 In this paper we will be using the class  $\Psi^{\nu}_{cl}(M)$ of classical operators (see e.g. \cite{Duis:BK-FIO-2011}).
 Furthermore, we denote by $\Psi_{+}^{\nu}(M)$ the class of positive definite operators in 
 $\Psi^{\nu}_{cl}(M)$,
 and by $\Psi_{e}^{\nu}(M)$ the class of elliptic operators in $\Psi^{\nu}_{cl}(M)$. Finally, 
 $\Psi_{+e}^{\nu}(M):=\Psi_{+}^{\nu}(M)\cap \Psi_{e}^{\nu}(M)$ 
 will denote the  class of classical positive elliptic 
 pseudo-differential operators of order $\nu$.
 
We associate a discrete Fourier analysis to the operator $E\in\Psi_{+e}^{\nu}(M)$ inspired by that
considered by Seeley (\cite{see:ex}, \cite{see:exp}), see also \cite{Greenfield-Wallach:hypo-TAMS-1973}. The eigenvalues of $E$ form a sequence $\{\lambda_j\}_{0\leq j<\infty}$, with $\lambda_0:=0 $ and with multiplicities $d_j$. The corresponding orthonormal basis of $L^2(M,dx)$ consisting of eigenfunctions of $E$ will be denoted by $\{e^k_j\}_{j\geq 0}^{1\leq k\leq d_j}$. Relative to this basis, Fourier coefficients, Plancherel identity and Fourier inversion formula can be obtained.

By introducing suitable Sobolev spaces on $M\times M$ adapted to $E$ we give first 
sufficient Sobolev type conditions on Schwartz integral kernels in order to ensure that 
the corresponding integral operators belong to
different Schatten classes. Second, we introduce notions of invariant operators, Fourier multipliers and full matrix-symbols relative to $E$. We apply those notions to characterise Schatten classes and to find sufficient conditions for nuclearity and $r$-nuclearity in the sense
of Grothendieck \cite{gro:me}
for invariant operators in terms of the full matrix-symbol.

 The problem of finding such criteria on different kinds of domain is classical and has been much studied (cf.  \cite{bs:sing}, \cite{dr13:schatten}, \cite{dr13a:nuclp}, \cite{Toft:Schatten-modulation-2008}, \cite{Buzano-Toft:Schatten-Weyl-JFA-2010} and
 references therein). 

\section{Kernel conditions for Schatten classes on compact manifolds}  
\label{SEC:Schatten-classes}
We first define Sobolev spaces $H^{\mu_1,\mu_2}_{x,y}(M\times M)$
of mixed regularity $\mu_1,\mu_2\geq 0$. 

\begin{defn}\label{def2} 
Let $K\in L^2(M\times M)$ and let $\mu_1,\mu_2\geq 0$.  
For
operators $E_j\in\Psi^{\nu_j}_{+e}(M)$
($j=1,2$) with $\nu_j>0$, we define
\begin{equation}\label{EQ:mixed-Sobolev}
K\in H^{\mu_1,\mu_2}_{x,y}(M\times M) \Longleftrightarrow
(I+E_1)_x^{\frac{\mu_1}{\nu_1}} (I+E_2)_y^{\frac{\mu_2}{\nu_2}} K \in L^2(M\times M),
\end{equation}
where the expression on the right hand side means that we are applying
pseudo-differential operators on $M$ separately in $x$ and $y$. We note that
these operators commute since they are acting on different sets of variables
of $K$. 
\end{defn}

By the elliptic regularity the spaces  $H^{\mu_1,\mu_2}_{x,y}(M\times M)$ do not depend
on a particular choice of $E_{1}, E_{2}$ as above.
We will now give our main kernel criterion obtained in \cite{dr:suffkernel} for Schatten classes. 
\begin{thm} \label{ext322} 
Let $M$ be a closed manifold of dimension $n$ and let $\mu_1, \mu_2 \geq 0$. 
Let 
$K\in L^2(M\times M)$ be such that $K\in H^{\mu_1,\mu_2}_{x,y}(M\times M)$. Then
the corresponding integral operator $T_K$ on 
 $L^2(M)$ 
 given by $T_{K}f(x)=\int_{M} K(x,y) f(y) dy$
 is in the Schatten classes $S_r(L^2(M))$ for $r>\frac{2n}{n+2(\mu_1+\mu_2)}$.
\end{thm}

We formulate below the results on the trace class. 
Due to possible singularities of the kernel $K(x,y)$ along the diagonal, 
we will require an averaging process on $K$ as described in \cite{dr:suffkernel}
(see also \cite{bri:k2}, \cite{del:tracetop}). 
We denote by $\widetilde{K}(x,y)$ the pointwise values of such 
process defined a.e. As a corollary of Theorem \ref{ext322}, for the trace class we have:
\begin{cor} \label{ext322aw} 
Let $M$ be a closed manifold of dimension $n$ and let $K\in L^2(M\times M)$, 
$\mu_1, \mu_2 \geq 0$, be such that  $\mu_1+\mu_2>\frac n2$ and 
$K\in H^{\mu_1,\mu_2}_{x,y}(M\times M)$. Then
the corresponding integral operator $T_K$ on 
 $L^2(M)$ is trace class and its trace is given by
\begin{equation}\label{EQ:trace}
\Tr(T_K)=\int\limits_M\widetilde{K}(x,x)dx.
\end{equation}
\end{cor}
We also obtain several corollaries in terms of the derivatives of the kernel. We denote 
by $C_x^\alpha C_{y}^{\beta}(M\times M)$ the space of functions of class $C^{\beta}$ 
with respect to $y$ and $C^\alpha$ with respect to $x$.

 \begin{cor} \label{cor25} 
 Let $M$ be a closed manifold of dimension $n$. Let   
 $K\in C_x^{\ell_1} C_{y}^{\ell_2}(M\times M)$ for some even integers $\ell_1,\ell_2\in 2\mathbb N_0$
 such that $\ell_1+\ell_2>\frac n2$. Then $T_K$ is in $S_1(L^2(M))$, and its trace is given by 
\begin{equation}\label{EQ:trace2}
\Tr(T_{K})=\int\limits_M K(x,x)dx.
\end{equation}
\end{cor}
The corollary above is sharp (cf. Remark 4.5 of \cite{dr:suffkernel}) as a consequence of classical results for the convergence of Fourier series on the torus (cf. \cite{ste-we:fa}, Ch. VII; \cite{wai:trig}).

We now formulate some consequences on compact Lie groups combining results from \cite{dr13:schatten}. For example, on the compact Lie group
$\SU2$ let $X,Y,Z$ be three left-invariant vector fields $X,Y,Z$ such that
$[X,Y]=Z$ (for example, these would be derivatives with respect to Euler angles at a point
extended to the whole of 
$\SU2$ by the left-invariance). Let $\mathcal L_{sub}=X^{2}+Y^{2}$ be the
sub-Laplacian. Then we have 

\begin{equation}\label{EQ:subL-Schatten-SU2}
0<r<\infty \textrm{ and } \alpha r>4 
\Longrightarrow (I-\mathcal L_{sub})^{-\alpha/2}\in S_r(L^{2}(\SU2)).
\end{equation}
The same is true for $\mathbb S^{3}\simeq \SU2$ considered as the compact
Lie group with the quaternionic product.
Using this instead of elliptic operators we can show:

\begin{cor} \label{COR:sub-Lap-SU2} 
Let $K\in L^2(\mathbb S^{3}\times \mathbb S^{3})$ be such that 
we have $(1-\subll)_x^{\mu_1/2}(1-\subll)_y^{\mu_2/2}K\in L^2(\mathbb S^{3}\times \mathbb S^{3})$ for some $\mu_1, \mu_2\geq 0$. 
Then $T_K$ is in  $S_r(L^2(\mathbb S^{3}))$ for $r>\frac{4}{2+\mu_1+\mu_2}$.
The same result holds on the compact Lie groups $\SU2$ and $\SO3$.
\end{cor}

We now argue that instead of the sub-Laplacian other globally hypoelliptic operators
can be used, also those that are not necessarily covered by H\"ormander's sum of
the squares theorem. We will formulate this for the group
$\SO3$ noting that, however, the same conclusion holds also on
$\SU2\simeq \mathbb S^{3}$. We fix three left-invariant vector fields $X, Y, Z$ on $\SO3$ 
associated to the derivatives with
respect to the Euler angles, 
so that we also have $[X,Y]=Z$,
see \cite{rt:book} or \cite{rt:groups} for the detailed expressions. We consider the following family of `Schr{\"o}dinger' differential operators
\[\mathcal H_{\gamma}=iZ-\gamma(X^2+Y^2),\]
for a parameter $0<\gamma<\infty.$ For $\gamma=1$ it 
was shown in \cite{Ruzhansky-Turunen-Wirth:arxiv} that
$\mathcal H_1+cI$ is globally hypoelliptic if and only if
$0\not\in\{c+\ell(\ell+1)-m(m+1): \ell\in\N, m\in\mathbb Z, |m|\leq\ell\}.$ 
%
%
It has been also shown in \cite[Section 4]{dr13:schatten} that, if $\gamma > 1$, then 
$I+\mathcal H_{\gamma}$ is globally hypoelliptic, and
$$(I+\mathcal H_{\gamma})^{-\alpha/2}\in S_p \textrm{ if and only if } \alpha p>4.$$
 
As a consequence of this and following the argument in \cite{dr:suffkernel} for the proof Theorem \ref{ext322} 
with $I+\mathcal H_{\gamma}$ instead of $E=\lapm$ for the manifold $M=\SO3$, as well
as Corollary \ref{COR:sub-Lap-SU2}, we obtain:

\begin{cor} \label{ext322so3}  
Let $K\in L^2(\SO3\times\SO3)$ be such that $(I+\mathcal H_{\gamma})_x ^{\mu_1/2}
(I+\mathcal H_{\gamma})_y ^{\mu_2/2}K\in L^2(\SO3\times\SO3)$ for some $\mu_1, \mu_2\geq 0$. 
Then the integral operator $T_{K}$ on 
 $L^2(\SO3)$ is in $S_r$ for $r>\frac{4}{2+\mu_1+\mu_2}$
 and $\gamma >1$. 
\end{cor}

\section{Symbols, Fourier multipliers and nuclearity}  
\label{SEC:Schatten-classes2}

Let us now consider the concepts of invariant operator and corresponding full symbols 
introduced in \cite{dr14a:fsymbsch}.  
The eigenvalues of $E\in\Psi_{+e}^{\nu}(M)$ (counted without multiplicities)
form a sequence $\{\lambda_j\}$ which we order so that
\begin{equation}\label{EQ:lambdas}
0=\lambda_{0}<\lambda_{1}<\lambda_{2}<\cdots.
\end{equation}
For each eigenvalue $\lambda_j$, there is
the corresponding finite dimensional eigenspace $H_j$ of functions on $M$, which are smooth due to the 
ellipticity of $E$. We set 
$$
d_j:=\dim H_j, 
\textrm{ and } H_0:=\ker E.
$$
We also set $d_{0}:=\dim H_{0}$. Since the operator $E$ is elliptic, it is Fredholm,
hence also $d_{0}<\infty$ (we can refer to \cite{Atiyah:global-aspects-1968}, \cite{ho:apde2} for
various properties of $H_{0}$ and $d_{0}$).
We fix  an orthonormal basis of $L^2(M)$ consisting of eigenfunctions of $E$:
\beq\label{fam}\{e^k_j\}_{j\geq 0}^{1\leq k\leq d_j},\eq 
where $\{e^k_j\}^{1\leq k\leq d_j}$ is an orthonormal basis of $H_j$. 
Let $P_j:L^2(M)\rightarrow H_j$ be the corresponding projection.

The Fourier coefficients of $f\in L^2(M)$ with respect to the orthonormal basis $\{e^k_j\}$ 
will be denoted by 
\begin{equation}\label{EQ:F-coeff}
(\efee f)(j,k):=\widehat{f}(j,k):=(f,e_j^k).
\end{equation}
We will call the collection of $\widehat{f}(j,k)$ the {\em Fourier coefficients of $f$ relative to $E$},
or simply the {\em Fourier coefficients of $f$}.
If $f\in L^2(M)$, we also write
\[\widehat{f}(j)=\left(\begin{array}{c}\widehat{f}(j,1)\\
\vdots\\
\widehat{f}(j,d_j)\end{array} \right)\in\ce^{d_j},\]
thus thinking of the Fourier transform always as a column vector.
The following theorem proved in \cite{dr14a:fsymbsch} is the base to introduce the concepts
 of invariant operators and full symbols relative to $E$.

\begin{thm}\label{THM:inv}
Let $M$ be a closed manifold and 
let $T:\cinfm\to L^{2}(M)$
be a linear operator.
Then the following
conditions are equivalent:
\begin{itemize}
\item[(i)] For each $j\in\ene_0$, we have $T(H_j)\subset H_j$. 
\item[(ii)]
For each $j\in\ene_{0}$ and $1\leq k\leq j$, we have
$TE e_{j}^{k}=ET e_{j}^{k}.$
\item[(iii)] For each $\ell\in\ene_0$ there exists a matrix 
$\sigma(\ell)\in\ce^{d_{\ell}\times d_{\ell}}$ such that for all $e_j^k$ 
\beq\label{invadef}\widehat{Te_j^k}(\ell,m)=\sigma(\ell)_{mk}\delta_{j\ell}.\eq
\item[(iv)]  For each $\ell\in\ene_0 $ there exists a matrix 
$\sigma(\ell)\in\ce^{d_{\ell}\times d_{\ell}}$ such that
 \[\widehat{Tf}(\ell)=\sigma(\ell)\widehat{f}(\ell)\]
 for all $f\in\cinfm.$
\end{itemize}
The matrices $\sigma(\ell)$ in {\rm (iii)} and {\rm (iv)} coincide.
If $T$
extends to a linear continuous operator
$T:\Dcal'(M)\rightarrow \Dcal'(M)$ then
the above properties are also equivalent to the following ones:
\begin{itemize}
\item[(v)] For each $j\in\ene_0$, we have
$TP_j=P_jT$ on $C^{\infty}(M)$.
\item[(vi)]  $TE=ET$ on $L^{2}(M)$.
\end{itemize}
\end{thm} 

If any of the equivalent conditions (i)--(iv) of Theorem \ref{THM:inv} are satisfied, 
we say that the operator $T:\cinfm\rightarrow L^{2}(M)$ 
is {\em invariant (or is a Fourier multiplier) relative to $E$}.
We can also say 
that $T$ is $E$-invariant or is an $E$-multiplier. When there is no risk of confusion we will just refer to such kind of operators 
as invariant operators or as multipliers.
If $T$
extends to a linear continuous operator
$T:\Dcal'(M)\rightarrow \Dcal'(M)$ then we will say that $T$ is 
{\em strongly invariant relative to $E$}.

The proposition below shows how invariant operators can be expressed in terms of their symbols.

\begin{prop}\label{PROP:quant}
An invariant operator $T_{\sigma}$ associated to the symbol $\sigma$ can be written 
in the following way:
$$
T_{\sigma}f(x)=\sum\limits_{\ell=0}^{\infty}\sum\limits_{m=1}^{d_{\ell}}(\sigma(\ell)\widehat{f}(\ell))_me_{\ell}^m(x)\label{form23}
=\sum\limits_{\ell=0}^{\infty}[\sigma(\ell)\widehat{f}(\ell)]^{\top} e_{\ell}(x),
$$
where $[\sigma(\ell)\widehat{f}(\ell)]$ denotes the column-vector, and 
$[\sigma(\ell)\widehat{f}(\ell)]^{\top}e_{\ell}(x)$ denotes the multiplication
(the scalar product)
of the column-vector $[\sigma(\ell)\widehat{f}(\ell)]$ with the column-vector
$e_{\ell}(x)=(e_{\ell}^{1}(x),\cdots, e_{\ell}^{m}(x))^{\top}$.
In particular, we also have
\begin{equation}\label{EQ:Tsigma-e}
(T_{\sigma}e_{j}^{k})(x)=\sum\limits_{m=1}^{d_{j}} \sigma(j)_{mk}e_{j}^{m}(x).
\end{equation}
If $\|\sigma(\ell)\|_{\mathcal L(H_{\ell})}$ grows polynomially in $\ell$
and $f\in C^{\infty}(M)$, the convergence in
\eqref{EQ:Tsigma-e} is uniform.
\end{prop}
We can now formulate our characterisation of the membership of invariant operators in Schatten classes:
\begin{thm}\label{schchr} Let $0<r<\infty$. 
An invariant operator $T:L^2(M)\rightarrow L^2(M)$ is in $S_r(L^2(M))$ 
if and only if its symbol $\sigma_{T}$ satisfies
$\sum\limits_{\ell=0}^{\infty}\|\sigma_{T}(\ell)\|_{S_r}^r<\infty$. 
Moreover
\[\|T\|_{S_r(L^2(M))}^r=\sum\limits_{\ell=0}^{\infty}\|\sigma_{T}(\ell)\|_{S_r}^r.\]
If an invariant operator $T:L^2(M)\rightarrow L^2(M)$ is in the trace class
$S_1(L^2(M))$, then 
\[\Tr(T)=\sum\limits_{\ell=0}^{\infty}\Tr(\sigma_{T}(\ell)).\]
\end{thm}


We now turn to some applications to the nuclearity on $L^p(M)$ spaces. 
Let $F_{1}$ and $F_{2}$ be two Banach spaces and $0<r\leq 1$, a linear operator $T$
from $F_{1}$ into $F_{2}$ is called {\em r-nuclear} if there exist sequences
$(x_{n}^\prime)\mbox{ in } F_{1}^{\prime} $ and $(y_n) \mbox{ in } F_{2}$ so that
\beq 
Ax= \sum\limits_n \left <x,x_{n}'\right>y_n \,\mbox{ and }\,
\sum\limits_n \|x_{n}'\|^{r}_{F_{1}^{\prime}}\|y_n\|^{r}_{F_{2}} < \infty.\label{rn}
\eq
This notion, developed by Grothendieck \cite{gro:me}, extends the notion of
Schatten classes to the setting of Banach spaces.

In order to study nuclearity on $L^p(M)$ for a given compact manifold $M$ of dimension $n$, we introduce a function $\Lambda(j,k;n,p)$ which controls the $L^p$-norms of the family of eigenfunctions $\{e_{j}^{k}\}$ of the operator $E$, i.e. we will suppose that 
$\Lambda(j,k;n,p)$ is such that we have the estimates
\beq\|e_{j}^{k}\|_{L^{p}(M)}\leq \Lambda(j,k;n,p).\label{control}
\eq
There are many things that can be said about the behaviour of $ \Lambda(j,k;n,p)$
in different settings, see e.g. results and discussions in
\cite{Donnelly:Asian-2006,Toth-Zelditch:DMJ-2002,dr14a:fsymbsch}.

We will use the following function 
$\tilde{p}$ for $1\leq p\leq \infty$: 
\beq\label{ppp} \tilde{p}:=\left\{
\begin{array}{rl}
0\,,&\, \mbox{ if } 1\leq p\leq 2,\\
\frac{p-2}{p},\,&\,\mbox{ if } 2<p<\infty ,\\
1,\,&\,\mbox{ if } p=\infty .
\end{array} \right. \eq
For $1\leq p_{1}, p_{2}\leq \infty$ we denote their dual indices by $q_{1}:=p_{1}^{\prime}$,
$q_{2}:=p_{2}^{\prime}$.
The criterion for $r$-nuclearity now is:
\begin{thm}\label{main126} 
Let $1\leq p_1,p_2 <\infty$ and $0<r\leq 1$. Let $T:L^{p_1}(M)\to L^{p_2}(M)$ be a
strongly invariant linear 
continuous operator. Assume that its matrix-valued symbol $\sigma(\ell)$ satisfies 
\[\sum\limits_{\ell=0}^{\infty} \sum\limits_{m,k=1}^{d_{\ell}}
|\sigma(\ell)_{mk}|^r\Lambda(\ell,m;n,\infty)^{\tilde{p_2}r}\Lambda(\ell,k;n,\infty)^{\tilde{q_1}r}<\infty.\]
Then the operator $T:L^{p_1}(M)\rightarrow L^{p_2}(M)$ is $r$-nuclear. 
\end{thm}

In some cases it is possible to simplify the sufficient condition above 
when the control function $\Lambda(\ell,m;n,\infty)$ is independent of $m$. 
For instance a classical result (local Weyl law) due to H{\"o}rmander 
(\cite[Theorem 5.1]{ho:eigen1},  \cite[Chapter XXIX]{ho:apde4}) implies the following
estimate:
\begin{lem} \label{LEM:Horm}
Let $M$ be a closed manifold of dimension $n$. Let $E\in \Psi_{+e}^{\nu}(M)$, then
\beq\|e_{\ell}^m\|_{L^{\infty}}\leq C\lambda_{\ell}^{\frac{n-1}{2\nu}}.\label{asein}\eq 
\end{lem}
Thus $\Lambda(\ell;n,\infty)=C\lambda_{\ell}^{\frac{n-1}{2\nu}}$ furnishes an
example of $\Lambda$ independent of $m$. For controls of type $\Lambda(\ell;n,\infty)$ we have a
basis-independent condition: 

\begin{cor}\label{main12ab} 
Let $1\leq p_1,p_2 <\infty$ and $0<r\leq 1$. 
Let $T:L^{p_1}(M)\to L^{p_2}(M)$ be a strongly invariant formally self-adjoint  
continuous operator. Assume that its matrix-valued symbol $\sigma(\ell)$ satisfies 
\[
\sum\limits_{\ell=0}^{\infty}
\|\sigma(\ell)\|_{S_r}^r\Lambda(\ell;n,\infty)^{(\tilde{p_2}+\tilde{q_1})r}<\infty.
\]
Then the operator $T:L^{p_1}(M)\rightarrow L^{p_2}(M)$ is $r$-nuclear. 
In particular, if its matrix-valued symbol $\sigma(\ell)$ satisfies 
\begin{equation}\label{EQ:cor-nuc}
\sum\limits_{\ell=0}^{\infty}
\|\sigma(\ell)\|_{S_r}^r\lambda_{\ell}^{{\frac{(n-1)}{2\nu}}(\tilde{p_2}+\tilde{q_1})r}<\infty,
\end{equation}
then the operator $T:L^{p_1}(M)\rightarrow L^{p_2}(M)$ is $r$-nuclear. 
\end{cor}
We now give an example of the application of such results
in the case of the sphere $\mathbb{S}^3\simeq \SU2$. 
We consider the Laplacian 
(the Casimir element) $E=-\mathcal{L}_{\mathbb{S}^3}$.

\begin{cor}\label{COR:su2-LapM}
If $\alpha > \frac 3r+\half(\tilde{p_2}+\tilde{q_1})$, $0<r\leq 1$, 
$1\leq p_1,p_2 <\infty$, the operator
$(I-\mathcal{L}_{\mathbb{S}^3})^{-\frac{\alpha}{2}}$ is $r$-nuclear from 
$L^{p_1}(\mathbb{S}^3)$ into $L^{p_2}(\mathbb{S}^3)$.
\end{cor}



\section*{Acknowledgements.}
The first author was supported by Marie Curie IIF 301599. 
The second author was supported by EPSRC grant EP/K039407/1.
The authors would like to thank V\'eronique Fischer for discussions.






\end{document}